\documentclass[fleqn, 11pt]{article}

\usepackage{fullpage} 
\usepackage[english]{babel}
\usepackage{amsmath, amssymb, amsfonts, theorem, dsfont}

\author{Misja Nuyens}
\title{The Foreground-Background  queue: a survey}
\newtheorem{theorem}{Theorem}[section]

\newtheorem{definition}[theorem]{Definition}

\newtheorem{corollary}[theorem]{Corollary}

\newcommand{\FB}{{\sf FB}}\newcommand{\FIFO}{{\sf FIFO}}\newcommand{\LIFO}{{\sf LIFO}}\newcommand{\PS}{{\sf PS}}\newcommand{\FBPS}{{\sf FBPS}}\newcommand{\SET}{{\sf SET}}\newcommand{\LAST}{{\sf LAST}}\newcommand{\LAS}{{\sf LAS}}\newcommand{\SRPT}{{\sf SRPT}}

\begin{document}
\maketitle
\begin{abstract}
The Foreground-Background (\FB) discipline, which gives service to the customer that has received the least amount of service, minimises the queue length for a certain class of heavy-tailed service times. In this paper we give an overview of the results in the literature on single-server queues with the \FB\ discipline.
\end{abstract}
{\bf Keywords:}  \FB, \FBPS, \LAS,  service discipline, queue length, sojourn time\\
{\bf AMS Subject Classification:} primary 60K25, secondary 90B22

\section{Introduction} 

Queues with heavy tailed distributions are of growing importance, as they seem to be good models for internet traffic, see for example Crovella and Bestavros \cite{bestavros} and Taqqu et al.\ \cite{taqqu}. 
For such queues,  First-in-first-out (\FIFO)  and other standard disciplines do not perform well, and instead one should consider service disciplines that are able to react properly on the arrival of very large jobs. The Shortest-Remaining-Processing-Time discipline (\SRPT) provides an improvement, but it uses the exact size of a job, and this knowledge is not always available. Among the disciplines that do not use knowledge of the job size,  the Foreground-Background (\FB) discipline is optimal in the following sense: for a certain class of heavy-tailed distributions, \FB\ stochastically minimises the queue length, see Theorem \ref{dfr} below. In this paper we give an overview of the literature on queues with the \FB\ discipline. 

The \FB\ discipline works according to the following priority rule: priority is given to the job  that has received the least amount of service. If there are $n$ such jobs, for some $n\in\mathbb{N}$, then they are served simultaneously, i.e., each of them is  served at rate $1/n$.       Alternatively formulated, if  the {\em age}  of a job is the amount of work it has received, then a server using the \FB\ discipline always serves the {\em youngest}  job(s).

It is perhaps surprising that until recently the \FB\ discipline has received  little attention in the literature: almost all results on the \FB\ queue date from after 1980. This lack of attention may have to do with the only relatively recent interest in queues with heavy-tailed characteristics, as well as with the difficulties that arise in the  analysis of  the \FB\ queue.     
The first important work on \FB\ queues was done by Schrage \cite{schr} and Kleinrock \cite{klein2} in the 1960s and 1970s. They mainly studied the sojourn time of a job of size $x$, and derived its mean and Laplace transform. Around 1980, Pechinkin \cite{pech}, Schassberger \cite{sch1} and Yashkov \cite{yas84} obtained expressions for the generating functional of the steady-state queue length. Yashkov \cite{yas92} provided a survey of most results known at the time. At the beginning of the 1990's, Righter and Shanthikumar \cite{rs, shasha} proved optimality results for the transient queue length. Since 2000, there has been a host of new results, for example, on the distribution of the sojourn time (Borst et al.\ \cite{borst} and Nuyens et al.\ \cite{nwz}), the slowdown (Harchol-Balter and Wierman \cite{mor2},  Rai et al. \cite{rai2003}), and the maximum queue length (Nuyens \cite{nuyens1}). 

In this paper  we intend to give an overview of the (different types of) theoretical results obtained for the single server \FB\ queue. A survey of this type does not exist in the literature. The survey articles \cite{yas87, yas92} by Yashkov on processor-sharing queues include results for \FB\ queues, but they date back from 1987 and 1992. Since then, both the interest in \FB\ queues and the number of results have strongly increased.  As a consequence, most results in the present survey are not contained in the older surveys. 
In addition to existing results, this paper contains some new material. Furthermore, to indicate similarities or differences, we shall compare the results for \FB\ with those for Processor-sharing (\PS), \SRPT\ and \FIFO.
We have chosen to focus on theoretical results for the M/G/1 and GI/GI/1 queue. Numerical results like those in Rai et al.\ \cite{rai2003} are not discussed, nor are the studies of multi-level queues with an \FB\ mechanism in Aalto et al.\ \cite{ayesta}. For readers interested in implementing \FB, we refer to Rai et al.\ \cite{rai2004, rai2005}.

This  survey  is organised as follows. Section \ref{anktj} is devoted to  the history of the \FB\ model and the acronym \FB, and a discussion of  general notions and intuition for the \FB\ queue. 
In Section \ref{ss:opt} we describe the optimality results for the \FB\ queue. The mean performance is considered in Section \ref{ss:meanperf}.
Section \ref{ss:ql} describes results on the stationary queue length, and the maximal queue length. Sojourn time asymptotics are the subject of Section \ref{sojjie}. The survey is concluded with Section \ref{slodo}, which is devoted to the slowdown.

\section{History of the \FB\ model  and used acronyms, Intuition and  general notions}\label{anktj}
Initially, in the second half of the 1960s, the term \FB, or rather \FB$_n$, was used as an abbreviation for both {\em Foreground-Background} and {\em Feedback} queueing systems. These different names referred to the same model, see Schrage \cite{schr}, Coffman and Kleinrock \cite{klein0}, and  the survey article on time-sharing models by McKinney \cite{mckinney}. The \FB$_n$ queue with so-called {\em quantum size} $q$ is a one-server queue with $n$ states, or priority classes. This queue operates as follows. 

Upon arrival in the queue, a job enters the first (or highest priority) state. Within each priority state, the priority of jobs depends on their arrival time to that state, in a \FIFO\ manner. Jobs are served one at a time  and uninterruptedly  for a time period of length $q$. After the server has completed a job's service request in a certain state, a job from the highest (non-empty) priority state is selected for service.
If a job does not leave the queue during its time in the $k$th 
state, it moves to state $k+1$ (which has lower priority) 
and waits until it is served in that state. In the $n$th and final state, jobs are served only if there are no jobs in other states. In that final state, they are served until they leave the system. So, for example, \FB$_1$=\FIFO.

The interest in the \FB$_n$ model with $n$ states and positive quantum size $q$ has faded a bit. Instead, people have studied  the limiting case where (first) $n\rightarrow \infty$ and (then) $q \to 0$. After  Kleinrock  devoted a section of \cite{klein2} to this limiting case of the \FB$_n$ model,  the term Foreground-Background (\FB) is generally used for this model, and so it is in this paper. We believe the term Foreground-Background is preferable over {\em Feedback}, since  there is no real feedback in the limit $q\to 0$.

To distinguish \FB\ from the \FB$_{n}$ model, some authors prefer to use the term Foreground-Background {\em Processor Sharing} (\FBPS), \FB$_{\infty}$
or {\em generalised} Foreground-Background. Others have invented their own acronyms like \LAS\ (Least Attained Service first) \cite{rai2003}, \LAST\ (Least Attained Service Time first) \cite{rs}, \cite{rsy}, or \SET\  (Shortest Elapsed Time) \cite{coffmandenning}.  Furthermore, \FB\ may be disguised as `advantageous sharing of a processor' \cite{pech} as well. Due to this cacophony of names,  some results on \FB\ queues are  difficult to find in the literature. 
One of the goals of the present survey is  to unite the \FB\ world again, and to provide a clear overview of all available results. This should prevent that theorems are  re-proved, sometimes even in a weaker form, like Theorem 2.1  in Feng and Misra  \cite{misra}.

\subsection{Intuition}
To get a feeling for the evolution of the queue under the \FB\ discipline, let us consider what happens  when a new job arrives to the  \FB\ queue.  Since  that job is (strictly) the youngest  in the queue,  it is served immediately. Now there are three scenarios. 
1. The new job needs at least as much service as the age of the job(s) that was (were) preempted at his arrival.  Then after some time  the job {\em joins} a {\em cohort}, a groups of jobs with the same age.  2. The job needs less service than the age of the job that was preempted.  Then the new job leaves the queue before joining the older cohort,  and the server returns to the cohort that was preempted.  3. Before joining another cohort or leaving the queue, the new job is preempted itself by the arrival of another new job.

Since a job with service time $x$ is younger than $x$ throughout its stay in the queue, it has {\em priority} over all jobs older than $x$. As a consequence, the time such a  job spends in the system  is  the same as if  all service times would be truncated at  $x$, i.e., if all service times $y$ would have value $\min\{y,x\}$ instead.  Hence, in the \FB\ queue,  small jobs do not suffer from the presence of large jobs in their midst. The sojourn times of small jobs are therefore {\em insensitive} to the shape of tail of the service-time distribution. This property turns out to be very useful  when studying sojourn times in  the \FB\ queue. 
It may be shown that the ratio of the expected time  a job spends in the system and its service time converges to 1 as the  service time converges to zero, see (\ref{kleinvx}) below. 

The consequence of this quick service to younger jobs is felt by jobs with (very) large demands. These jobs are mostly served when no other jobs are present,  see also the discussion below equality (\ref{corvx}).
One of the main issues when studying  \FB\ queues is to determine the price that large jobs have to pay for the priority that  \FB\ gives to short jobs.

\subsection{Set-up}
Unless stated otherwise, we consider 
an M/G/1 queue with  arrival rate $\lambda$, service-time distribution $F$, and generic service time $B$. Assume that  the load $\rho$ satisfies $\rho=\lambda EB<1$. If $F$ has a density, we denote it by $f$. Let $Q$ denote the stationary queue length, measured in number of jobs. The busy period length is denoted by $L$. An important role is played by queues with truncated service times $B\wedge x=\min\{B,x\}$. The load $\rho(x)$ of such a queue is given by $\rho(x)=\lambda E(B\wedge x)$.
The class  ${\cal D}$ is the class of all service disciplines that do not use knowledge of the residual service times. So, for example, $\FB, \PS, \FIFO\in {\cal D}$, but $\SRPT\notin {\cal D}$.

Finally, some results in  this paper are stated in terms of certain  classes of service-time distributions, like DFR, and certain stochastic orderings. Readers unfamiliar with these notions are referred to the appendix for a discussion of these notions.

\section{Optimality of \FB}\label{ss:opt}
In this section we describe two optimality results for the transient queue length. Originally, they were proved in the classical GI/GI/1 setting. However, they hold true in the following more general set-up as well.

Consider queues with the same arrival and service times, but with different service disciplines.  Let $Q_{\pi}(t)$ denote the queue length at time $t$ in the queue with discipline $\pi$. The sequence of arrival times may be any  sequence, even deterministic. Furthermore, at time 0 some jobs could be present, as long as the number of jobs present at time 0 in  the two queues is the same and their ages at time 0 are (pairwise) equal. Finally, jobs may have received service prior to their arrival in the queue, i.e., they may enter the queue with a positive age.

The first optimality result, proved in Corollaries 2.1.2 and 2.3.2 in Righter and Shanthikumar \cite{rs}, is in terms of the marginal distributions of the process $\{Q(t), t\geq 0\}$. It generalises the result mentioned in Yashkov \cite{yas87} that for DFR distributions, \FB\ minimises $EQ$ over all disciplines in ${\cal D}$. 

\begin{theorem}\label{dfr} Consider a GI/GI/1 queue. Let $\pi\in{\cal D}$, and let $\sigma\in {\cal D}$ be a non-preemptive discipline.
If the service-time distribution belongs to the class DFR, then for  every $t\geq 0$,
\[ Q_{\FB}(t)\leq_{st}  Q_{\pi}(t)\leq Q_{\sigma}(t).\]
For IFR service times the inequalities are reversed.
\end{theorem}

Service-time distributions are often divided in heavy-tailed and light-tailed distributions, and so it is in this paper, for example in Section \ref{sojjie} below. Theorem \ref{dfr} indicates that for the \FB\ queue an alternative approach could be used, namely by considering  DFR and IFR distributions, although these two classes do not contain all possible distributions. As remarked in Appendix \ref{fiokje} below, there is a connection between these two divisions: some well-known heavy-tailed distributions are DFR, while some well-known light-tailed distributions are IFR.\\ 


Using a stronger condition on the density of the service-time distribution, an optimality result for the law of the process  $\{Q(t), t\geq 0\}$ can be obtained.  Theorem 13.D.8 by Righter  in \cite{rs}   reads:

\begin{theorem}\label{simthm} Consider a GI/GI/1 queue.
Let $\pi\in{\cal D}$, and let  $\sigma\in {\cal D}$ be a non-preemptive discipline.  If the service-time distribution has a log-convex density, then 
\begin{equation}\label{dlreq} \{Q_{\FB}(t), t\geq 0\} \leq_{st} \{Q_{\pi}(t) , t\geq 0\}\leq_{st}  \{Q_{\sigma}(t), t\geq 0\}.\end{equation}
For service times with a log-concave density, the inequalities are reversed.
\end{theorem}

Let us conclude with a few remarks on the proofs of Theorems \ref{dfr} and \ref{simthm}.  Theorem \ref{simthm} essentially already appeared in Righter and Shanthikumar  \cite{rsfifo}, but there it  was formulated  only in terms of the \FIFO\ discipline. The proofs of Theorems \ref{dfr} and \ref{simthm} are given in discrete time only. The extension to continuous time  is said to follow from a limit argument.  However, the limit argument given in chapter 10 of \cite{nuyenstheziz} suggests that such an argument is far from trivial.

A  third optimality result, Theorem 3.14 in Righter et al.\ \cite{rsy} states that if $E[B-x|B>x]$ is increasing in $x$,  then $ EQ_{\FB}(t)\leq  EQ_{\pi}(t)$ for all $t\geq 0$ and  $\pi\in {\cal D}$.  However, the proof contains an error that cannot be immediately fixed, as was noted by  Aalto et al.\  \cite{ayesta}:  the result is proved by considering the unfinished work of jobs with age less than $x$, but the proof does not take into account that this quantity makes a vertical downward jump  when a job reaches age $x$. The proof of  Lemma 2.4 in Feng and Misra \cite{misra} contains the same mistake.

Finally, the \FB\ discipline is in a certain  sense opposite to non-preemptive disciplines (i,e., disciplines that do not interrupt the service of a customer once is has started) such as \FIFO\ and \LIFO:  \FB\ serves the youngest jobs, while non-preemptive disciplines  give priority to the oldest job.  This idea is illustrated by the  optimality theorems.

\section{Mean performance}\label{ss:meanperf}
In this section we give several results that describe the mean performance of the \FB\ queue.  We consider the  mean queue length in the stationary queue, its heavy-traffic behaviour,  the influence of variability in the service-time distribution, and the conditional expected sojourn time.

Many of the results described in this survey illustrate the idea that for heavy-tailed distributions the \FB\ discipline performs very efficiently. Many results indicate that for those distributions  the \FB\ discipline is markedly superior to \FIFO. Furthermore, we shall encounter the following  perhaps counter-intuitive phenomenon: for heavy-tailed service times, the \FB\ queue may behave better than for light-tailed service times (see e.g.\ Corollary \ref{jarrow} below).

\subsection{Mean queue length}
The mean queue length in the M/G/1 \FB\ queue is given by the following expression, see
Pechinkin \cite{pech} - although there  the factor $\lambda$ in front of the integral is missing:

\begin{corollary} In the M/G/1 \FB\ queue,
\begin{equation}\label{sojtijd}EQ=\lambda\int_0^{\infty}\Big( \frac{\lambda E(B\wedge x)^2} 
{2(1-\rho(x))^2}+\frac{x}{1-\rho(x)}\Big) dF(x).\end{equation}
\end{corollary}
The same expression can be found by using Little's law in combination with Theorem \ref{samefash} below. Using that $EB<\infty$, it may be seen  from  (\ref{sojtijd})  that $EQ<\infty$.   In fact, Yashkov \cite{yas92} gives the following equality for the  expected queue length:
\begin{equation}\label{danasta}EQ\leq \frac{\rho (2-\rho)}{2(1-\rho)^2}.\end{equation}
In (\ref{danasta}),  equality holds if and only if the service-times are deterministic. It is not surprising that deterministic service times maximise the mean queue length: remember that in the M/D/1 \FB\ queue, all jobs leave at the end of the busy period.
For light-tailed service times, the behaviour of  the queue under the \FB\ discipline may be  quite poor. Consider, for example, a busy period in the G/D/1 queue. By the \FB\ priority rule,  none of the jobs is allowed to leave the queue before any other job.  Hence, all jobs that arrive during  a  busy period leave the queue together at the end of the busy period. Kleinrock \cite{klein2}  uses this example to emphasise the  disastrous effect that using the \FB\ discipline may have.  

Next we consider higher order moments of the stationary queue length. In the M/G/1 \PS\  queue, all moments of the stationary queue length $Q_{\PS}$ exist, since $P(Q_{\PS}=n)=(1-\rho)\, \rho^n$ for $ n\geq 0$. Furthermore, by  using the Pollaczek-Khinchin transform for the queue length, it may be seen that in the M/G/1 \FIFO\ queue,  $EQ_{\FIFO}^n<\infty$ if and only if $EB^{n+1}<\infty$. For higher order moments in the \FB\ queue, Theorem 9.1 in \cite{nuyenstheziz} gives the following relation between the moments of $B$ and $Q$.

\begin{theorem}\label{mommie}  In the M/G/1 \FB\ queue, if  $EB^{\alpha}<\infty$ for some $\alpha>1$,  then all moments of $Q$ are finite. 
\end{theorem}

It is an open question whether  $EB<\infty$ is enough to show that all moments of $Q$ exist. \\

\subsection{Heavy traffic}\label{ss:heavy}
In this subsection we describe  results on  the heavy-traffic behaviour of the mean queue length in the stationary  queue. Most of these results have the following form: for a certain  class of service-time distributions, there is a constant $\gamma>0$ such that $EQ=O((1-\rho)^{-\gamma})$ as $\rho\uparrow 1$. 

Let $x_F$ be the right end-point of the service-time distribution $F$, i.e., $x_F=\sup\{x: F(x)<1\}.$ It turns out that the mean queue length  shows different  heavy-traffic behaviour for $x_F=\infty$ and $x_F<\infty$. We start with the case that $x_F=\infty$. Recall that a function  $\eta$ is called  regularly varying at $\infty$ of index $\alpha$ if  $\lim_{x\to\infty} \eta(tx)/\eta(x)= t^{\alpha}$ for all $t>0$.

\begin{theorem} \label{ank2} Let $x_F=\infty$.
Then $(1-\rho)EQ=O(1)$ for $\rho\uparrow 1$ if one  of the following conditions holds: 
\begin{enumerate}
\item  $1-F$ is regularly varying at $\infty$ of index $\alpha<-1$,
\item $1-F(x)=o(\exp(-cx^{\beta}))$ for some $\beta>0$ as $x\to\infty$.
\end{enumerate}
\end{theorem}
This theorem was proved in \cite{nuyenstheziz}.
By comparing  the \FB\ queue with the \PS\ queue, and using Theorem \ref{dfr},  for DFR distributions we have $(1-\rho)EQ=O(1)$ holds as well.

Bansal and Gamarnik \cite{bansalgamarnik} obtained the following stronger results for Pareto distributions.
\begin{theorem}\label{bansalthm}  Let $1-F(x)=(k/x)^{\alpha}, x\geq k$, for some $\alpha>1$.  Then for $\rho\uparrow 1$, 
\[ EQ=\begin{cases}
O(\log(\frac{1}{1-\rho})) & {\rm if} \ \alpha<2,\\
O(\log^2(\frac{1}{1-\rho})) & {\rm if} \ \alpha=2,\\
O((1-\rho)^{-\frac{\alpha-2}{\alpha-1}}) & {\rm if} \ \alpha>2.
\end{cases}\]
\end{theorem}
By partial integration of  (\ref{sojtijd}), Yashkov \cite{yas92} showed that
$EQ\geq -\log (1-\rho)$ for $\rho<1$. Combining this with Theorem \ref{bansalthm}, we find the following new result.
\begin{corollary}  Let $1-F(x)=(k/x)^{\alpha}, x\geq k$, for some $1<\alpha<2$. 
Then there exist  $c_2\geq c_1\geq 1$ such that 
\[(c_1+o(1))\log \Big(\frac{1}{1-\rho}\Big)\leq EQ\leq (c_2+o(1))\log \Big(\frac{1}{1-\rho}\Big),\qquad \rho\uparrow 1. \] 
\end{corollary}

In addition,  a few other heavy-traffic results exist for the \FB\ queue. The survey paper \cite{yas92} quotes the following heavy-traffic  limits from articles that appeared in Russian journals. Unfortunately, these could  not be retrieved by the author of this survey. 
The following theorem is by Nagorenko and Pechinkin \cite{nago}.
Recall that $f(x)\sim g(x)$ for $x\to\infty$ means that $\lim_{x\to\infty} f(x)/g(x)=1$.
\begin{theorem} For service-time distributions with tail $1-F(x)\sim a x^b e^{-cx},$ for some $a>0, b\geq 0$ and $c>0$, the stationary queue length $Q$  in the M/G/1 \FB\ queue satisfies
\[ \lim_{\rho\uparrow 1} P(Q/EQ<x)=1-e^{-x}, \qquad x\geq 0.\]
\end{theorem}
For the M/D/1 \FB\ queue  an expression for  the Laplace transform of $\lim_{\rho\uparrow 1} Q/EQ$ is given by Yashkov and Yashkova \cite{yasenyas}.\\

We now turn to the case that $x_F<\infty$. It turns out that the heavy-traffic behaviour of the mean queue length is different than for the case that $x_F=\infty$.
\begin{theorem}\label{giko} Let $x_F<\infty$. If $1-F(x)\sim (x_F-x)^{\beta}\eta$   for some constants $\beta, \eta>0$ as $x\to x_F$,  then for $\rho\uparrow 1$,
\[(1-\rho)^{(\beta+2)/(\beta+1)}EQ= O(1).\]
\end{theorem}

Furthermore, for the M/D/1 queue we have the following result. By combining  the Pollaczek-Khinchin mean value formula and (\ref{danasta}), it can be seen  that the queue lengths $Q_{\FIFO}$ and $Q_{\FB}$ in the M/D/1 queue under \FB\ and \FIFO\ satisfy:
\[\frac{EQ_{\FIFO}}{EQ_{\FB}}=\frac{\rho}{2(1-\rho)}\frac{2(1-\rho)^2}{2-\rho}=
\frac{\rho(1-\rho)}{2-\rho}=1-\rho+O((1-\rho)^2), \qquad \rho\uparrow 1.\]

\subsection{Variability of the  service-time distribution}\label{compar}
In this section we consider the effect on the queue length of more variability in the service times.  
For heavy-tailed service times the \FB\ disciplines performs quite well. One may therefore wonder  whether in the \FB\ queue more variability in the service times could be beneficial to the behaviour of the queue; could it, for example,  lead to a smaller mean queue length?
In the literature, a few conjectures of this type exist.
 The survey paper Yashkov \cite{yas92} claims that in the stationary \FB\ queue   ``$EV$ decreases with an increase in the dispersion of $F(x)$, and conversely increases as the dispersion of $F(x)$ decreases.''

Furthermore, for or a random variable  $X$  with $EX> 0$, the {\em coefficient of variation} is defined as $C(X)=\sqrt{\textrm{Var}(X)}/EX$. Coffman and Denning \cite{coffmandenning} conjectured that an \FB\ queue with  generic service-time $B$ and $C(B)> 1$ has a smaller expected sojourn time than a queue  with generic service time $B'$ with the same mean and $C(B')<1$. Wierman et al.\  \cite{wierman1} invalidated this conjecture. 

However, a result somewhat similar to the conjecture does hold.
First note that by inequality (\ref{danasta}), for fixed $\rho$, the mean queue length  in the \FB\ queue  is maximal for deterministic service times. These have coefficient of variation zero. By using that in the M/M/1 queue all disciplines have the same queue-length distribution, we obtain the following corollary to Theorem \ref{dfr}.

\begin{corollary}\label{jarrow}
Consider four stationary M/G/1 \FB\ queues with arrival rate $\lambda$ and load $\rho$.  The service-time distributions are DFR, exponential, IFR, and  deterministic, respectively. Let $Q^{DFR}$, $Q^{exp}$, $Q^{IFR}$ and $Q^{det}$ denote the stationary queue lengths. Then
\[EQ^{DFR}\leq EQ^{exp}=\frac{\rho}{1-\rho}\leq EQ^{IFR}\leq 
EQ^{det}=\frac{(2-\rho)\rho}{2(1-\rho)^2}.\]
\end{corollary}
It can be shown that if $X\in$ DFR, $Y\in$ IFR, $Z$ is deterministic, and $EX=EY=EZ$, then  $C(X)\geq C(Y) \geq C(Z)$. Hence, for special classes of service-time distributions, the queue with the larger coefficient of variation does have the smaller mean queue length.  \\


\subsection{Conditional mean sojourn time}\label{mien}
Let $V(x)$ be the (conditional) sojourn time of a job with service time $x$ in the stationary M/G/1 \FB\ queue.   By analysing the mean behaviour of the queue, Schrage \cite{schr} found the following expression for  $EV(x)$.
\begin{theorem}\label{samefash}
For all $x$ such that $\rho(x)<1$, 
the conditional sojourn time $V(x)$ satisfies
\begin{equation}\label{kleinvx}
EV(x)=\frac{\lambda E(B\wedge x)^2}{2(1-\rho(x))^2}+\frac{x}{1-\rho(x)}.
\end{equation}
\end{theorem}
A formula for $EV$ can be derived by integrating  this expression over all $x$ w.r.t.~the service-time distribution, see also (\ref{sojtijd}).   By differentiating (\ref{kleinvx}), Kleinrock  \cite{klein2}, obtained  the following consequence of Theorem~\ref{samefash}, under the condition $EB^2<\infty$.  This condition was removed in \cite{nuyenstheziz} by using that $E(B\wedge x)^2=o(x)$ and $1-F(x)=o(1/x)$ if  $EB<\infty$.
\begin{corollary} The mean conditional sojourn time $EV(x)$ satisfies
\begin{equation}\label{corvx} \lim_{x\rightarrow \infty} \frac{dEV(x)}{dx}=\frac{1}{1-\rho}.\end{equation}
\end{corollary}
This result may be interpreted as follows: during the sojourn time of an exceptionally large job, the service rate it gets is the total service rate (namely 1) reduced by the load of jobs that pass through the system in the meantime ($\rho$, in the limit).

Combining Theorem 1 and Lemma 1 in Rai et al.\  \cite{rai2003} yields  the following result.
\begin{theorem} For all $x$, 
\[ EV(x)\leq \frac{(2-\rho)x}{2(1-\rho)^2},\]
and equality holds if and only if $P(B=x)=1$.
\end{theorem}

Furthermore, Kleinrock \cite{klein2} gives an expression for the Laplace transform of $V(x)$:
\begin{equation}\label{kleinwait} E\exp(-sV(x))=\exp\Big(-x(s+\lambda-\lambda g_x(s))\Big) E\exp\Big(-(s+\lambda-\lambda g_x(s))W(x)\Big),\end{equation}
where $W(x)$ is the stationary workload in the queue with generic service time $B\wedge x$, and $g_x$ is the Laplace transform of the busy period in such a queue. 
By differentiation of these Laplace transforms,  in \cite{nuyenstheziz}  the following asymptotics  for $x\to\infty$ are derived: 
\begin{equation}\label{hopla} EV(x)^n=\frac{x^n}{(1-\rho(x))^n}+\begin{cases} O(x^{n-1}) &  {\rm if} \ \exists \, \alpha \geq 2 : EB^{\alpha}<\infty,\\
 o(x^{n+1-\alpha}) & {\rm if} \ \exists \, 1<\alpha<2 : EB^{\alpha}<\infty.\end{cases}\end{equation}
Since trivially $V(x)\geq x$, using (\ref{hopla})  yields the following new result.
\begin{theorem}\label{eveb}
For any $n\in\mathbb{N}$, we have  $EV^n<\infty$ if and only if $EB^n<\infty$.
\end{theorem}
This theorem indicates that the tail behaviour of the sojourn time and the service time is in some sense similar. This observation is further illustrated in Theorem \ref{sinthm1} below.

\subsection{Overload}
Assume $\rho>1$. Then the workload in a queue grows a.s.~without bounds. Under some service disciplines, e.g., \FIFO, still  all jobs will  eventually leave the system. In the \FB\ queue, only jobs with  service time below  a critical value a.s.\ leave the system. Indeed,  jobs smaller than  $x^*$, with
\[ x^*=\sup\{x: \rho(x)=\lambda E(B\wedge x)<1\},\] 
experience a system for which the stability condition $\rho(x)<1$ holds, and hence they leave the queue a.s. Jobs with service time larger than the critical value have a positive probability of being in the queue forever.

For the \PS\ queue, an expression for the asymptotical growth rate is given by Jean-Marie and Robert \cite{robert}.  It is interesting to compare these two growth rates. 
Numerical calculations indicate that for some  heavy-tailed distributions with much mass to the left of $\lambda^{-1}$, the asymptotic growth rate under \FB\ is smaller than under \PS.  However, no  theoretical results have been obtained so far.

\section{The queue length}\label{ss:ql}
In this section we describe the remaining results for the queue length.
In Section \ref{ss:max} the maximal queue length is treated; in Section \ref{ss:stat} we give transforms of the stationary queue length, and results for the cohort process.
\subsection{The maximal queue length} \label{ss:max}
Now we consider the maximal queue length in a busy period, $M$. The distribution of $M$ in the M/D/1 queue essentially already appears in Borel \cite{borel}:

\begin{theorem}\label{prabbie} 
The distribution of the maximal length $M$ in the  M/D/1 \FB\ busy period with  arrival rate $\lambda$ and service times equal to 1 satisfies
\[P(M=n)=\lambda^{n-1}e^{-\lambda n} \frac{n^{n-1}}{n!}\sim e^{n(\log \lambda +1-\lambda)}/(n\sqrt{n}\lambda\sqrt{2\pi}), \quad n\to\infty.\]
\end{theorem}
Here $f(n)\sim  g(n)$ means that $\lim_{n\to\infty}f(n)/g(n)=1$.

Nuyens \cite{nuyens1} proved that if the service times have a log-convex density, then the tail of  $M$  is bounded by an exponential tail:
\begin{theorem}\label{rsktj}
 If the service-time distribution has a log-convex density, then 
\[P(M>n)\leq \rho^n, \qquad  n=0,1,\ldots.\] \end{theorem}

Interestingly,  given that the service-time distribution has a log-convex density, the upper bound $\rho^n$ for $P(M>n)$ is insensitive to the precise form of the service-time distribution. Furthermore, one may wonder whether Theorem \ref{rsktj} implies that  all moments of $Q$ are finite for log-convex densities, see also Theorem \ref{mommie}.\\

By the regenerative structure of the queue length process,
the maximal queue length during a busy period is related to the maximal queue length over the time interval $[0,t]$ for $t\to\infty$, see the survey article Asmussen \cite{asmussen}. Nuyens \cite{nuyens1} used this idea to show the following.
\begin{theorem}\label{love}
Let $M(t)$ be the maximal queue length over the time interval $[0,t]$  of an M/G/1 \FB\ queue with i.i.d.~service times with a log-convex density. Assume that $\rho<1$. Then for any $x> 0$, the inequality
\[ P(M(t)>a \log t+b+x)\leq \rho^x,\]
holds for $t$ large enough, where $a=-1/(\log \rho)$, $b=-(\log \lambda+\log (1-\rho))/(\log\rho)+1$.
\end{theorem}
Using Theorem \ref{love}, calculations in \cite{nuyens1} showed that for some heavy-tailed log-convex densities, the time to overflow of a buffer in the \FB\ queue is of a different  order of magnitude than in the \FIFO\ queue. This illustrates the idea that in case of heavy-tailed service times, using the \FB\ discipline instead of \FIFO\ may increase the performance of the queue considerably.

\subsection{The stationary queue length}\label{ss:stat}
Pechinkin  \cite{pech} obtained the following expression for the  generating function of $Q$. 
\begin{theorem}\label{ptsj} Let $Q$ be the  number of jobs in the
stationary  M/G/1 queue. Then for  $z<1$,
\begin{equation}\label{ezn} Ez^Q=(1-\rho)\exp\Big({-\int_0^{\infty}z\frac{\partial v(t,z)}{\partial z}dt}\Big),\end{equation}
where $v(t,z)$ is the unique non-negative root of the equation
\begin{equation}\label{vtz} v(t,z)=\lambda \big(1-\int_0^te^{-v(t,z)x}dF(x)-z(1-F(t))e^{-v(t,z)t}
\big).\end{equation}
\end{theorem}
Yashkov \cite{yas84} obtained the counterpart of (\ref{ezn}) for the case of batch arrivals. From the proof of Theorem \ref{ptsj} it follows that  $v(t,1)=0$ and that $v$ is differentiable w.r.t.\ $z$. This allows for computing the moments of $Q$ by differentiating (\ref{ezn}). 

Let $Q_x$ be the number of jobs younger than $x$ in the stationary queue.
Schassberger  \cite{sch1} obtained the generating functional of the point process $Q_x$.

\begin{theorem} For all  functions $h:[0,\infty)\to [0,\infty)$,  \begin{equation}\label{stem} E(\exp(-\int h(x)dQ_x))=(1-\rho)\exp\big(-\int_0^{\infty} z\frac{\partial v(t,z)}{\partial z}\Big|_{z=\exp(-h(t))}dt\big),\end{equation}  where $v(t,z)$ is again the unique nonnegative root of (\ref{vtz}). \end{theorem}

The original proof of (\ref{stem}) in  Schassberger  \cite{sch1} uses a discrete approximation and is quite technical. Later, Robert and Schassberger \cite{robsch} found a more direct way to prove (\ref{stem}), using (\ref{ezn}) in combination with the following result.\begin{theorem}\label{indepprop}
The process $Q_x$ has independent increments.
\end{theorem}
The process  $Q_x$ is integer-valued, but it  has a jump  larger than 1 at $x$ when the cohort of age $x$ contains more than one job. Hence, it  is a {\em non-homogeneous} Poisson process with {\em batch arrivals}, see  Daley and Vere-Jones \cite{daleyvere}.  
Let  $K_x$ be the number of {\em cohorts} in the stationary M/G/1 \FB\ queue consisting of jobs younger than  $x$. The cohort process $K_x$ has only jumps of size 1.  Combining this observation with  Theorem \ref{indepprop}  leads to the following new corollary.
\begin{corollary} The cohort process  $K_x$ is a non-homogeneous Poisson process with intensity $\mu(x)$ given by
\[ \mu(x)=-\frac{d}{dx} \log P(K(x)=0)=-\frac{d}{dx} \log P(Q_x=0)=-\frac{d}{dx} \log (1-\rho(x))=\frac{\lambda (1-F(x))}{1-\rho(x)}.\]
\end{corollary}

\section{The sojourn time}\label{sojjie}
In this section we describe the  results obtained for the tail behaviour of the sojourn-time distribution.  In Section \ref{teq} we discuss the case of heavy-tailed sojourn times, and in Section \ref{lteq} we  consider light-tailed service times.

\subsection{Asymptotics: tail equivalence for heavy tails}\label{teq} 
For a broad class of heavy-tailed service times, N\'{u}\~{n}ez Queija  \cite{sindothes} obtained asymptotics for the tail of the service time distribution. Recently, Nuyens et al.\ \cite{nwz} generalised this result to the GI/GI/1 setting. A key role is played by the following class of functions. A function $h$  is said to be of {\em intermediate regular variation at infinity}, if 
\[\liminf_{\varepsilon\downarrow 0} \liminf_{x\rightarrow \infty}
\frac{h(x(1+\varepsilon))}{h(x)}=1.\]

\begin{theorem}\label{sinthm1}
Suppose that $1-F$  is of intermediate regular variation at infinity. If   $EB^{\alpha}<\infty$ for some $\alpha>1$,  then
\begin{equation}\label{eqsinthm1}  P(V >x)\sim P(B>(1-\rho)x).\end{equation} \end{theorem}
Theorem \ref{sinthm1}  makes sense intuitively: since jobs under \FB\ are only served if no younger jobs are present, very old jobs most likely only get service if no other jobs are present. Hence,  long jobs are served as if they were alone in a system with service rate $1-\rho$. This phenomenon is often called the {\em reduced load approximation.}  See also the remark after equation (\ref{corvx}). Furthermore, Theorem \ref{sinthm1} supports the remark below Theorem \ref{eveb}.

N\'{u}\~{n}ez Queija \cite{sindothes} showed that Theorem \ref{sinthm1} holds for M/G/1 \PS\ and \SRPT\  as well. Guillemin et al.\ \cite{grz} showed that relation (\ref{eqsinthm1}) holds for a wide class of processor sharing queues. For an overview the tail behaviour of $V$ under several other disciplines in case of regularly varying service times, see Borst et al.\ \cite{borst}. 

 It is an interesting open problem whether there exists a service discipline $\pi$ and some constant $1-\rho<a\leq 1$  such that in the  M/G/1 queue with heavy-tailed service times,
\begin{equation}\label{stretch}P(V_{\pi}>x)\sim P(B>ax).\end{equation}
If such a discipline does not exist, then the \FB\ queue has the perhaps counterintuitive  quality that although  large jobs are discriminated, large sojourn times are as unlikely as possible.
Another interesting question is how (\ref{stretch}) relates to the results in  Balter et al.\ \cite{mor1}, who study the behaviour of $V(x)/x$ for large $x$. These results are discussed in Section \ref{slodo} below.\\

Finally, let us compare the tail behaviour of the sojourn time under \FB\ and the busy period in the M/G/1 queue. Let $\eta$ be a regularly varying function at $\infty$ of index $\alpha<-1$. Let $L$ denote the busy period length. De Meyer and Teugels \cite{demeyer} proved that, as $x\to\infty$,
\begin{equation}\label{sonku} P(B>x)\sim \eta(x)\ \ \Leftrightarrow \ \ P(L>x)=(1-\rho)^{-\alpha-1} \eta(x).\end{equation}
Since $\eta(x)$ is also of intermediate regular variation at $\infty$, we have by Theorem \ref{sinthm1} and (\ref{sonku}),
\begin{align*} P(V>x)\sim P(B>(1-\rho)x) \sim \eta(x)((1-\rho))^{-\alpha}\sim (1-\rho) P(L>x).\end{align*}
Hence, in this case the tails of $V$ and $L$ are asymptotically equal up to a multiplicative constant. This property is shared with \PS\ and \SRPT, but it does not hold for $\FIFO.$

\subsection{Asymptotics: light-tailed service times}\label{lteq}
For light-tailed  distributions,  the behaviour of the sojourn-time distribution is  formulated in terms of the following quantity.
 \begin{definition}[Decay rate]\label{drdef} The {\em (asymptotic) decay rate} of a random variable $X$ is defined as $\gamma(X)=-\lim_{x\to\infty} x^{-1}\log P(X>x)$, provided the limit exists.\end{definition}
Mandjes and Nuyens \cite{mandjesnuyens} studied the decay rate of the sojourn time in the M/G/1 queue. Nuyens et al.\ \cite{nwz}  obtained the following generalisation.
\begin{theorem} \label{tiejorum}
Consider a GI/GI/1 \FB\ queue.
Let $L$ be the length of a busy period. Assume that $E(\exp(\kappa B))<\infty$ for some $\kappa>0$. Then  $\gamma(V)$  exists and
\[\gamma(V)= \gamma(L)=\sup_{s\geq 0}\{s+\Phi_A^{-1}\big(\frac{1}{\Phi_B(s)}\Big)\},\]
where $\Phi_B$ and  $\Phi_A$ are the generating functions of $B$ and the generic interarrival time $A$, and $\Phi_A^{-1}$ is the  inverse of $\Phi_A$. 
\end{theorem}

From the proof of Theorem 4.1 in Nuyens and Zwart \cite{nuyenszwart}, it follows that under any work-conserving discipline, i.e.\ a discipline that serves at rate 1 as long as there is work in the queue, $\gamma(V)\geq \gamma(L)$, if  the service times with an exponential moment. Hence, 
for service times with an exponential moment, the \FB\ discipline minimises the decay rate of the sojourn time in the class of work-conserving disciplines.

Contrasting Theorems \ref{dfr} and \ref{tiejorum}, we observe the following interesting phenomenon: for gamma-distributed service times  with a log-convex density, see Appendix  \ref{fiokje},   the \FB\ discipline minimises the queue length, but the sojourn time has the smallest possible decay rate. Hence, optimising one characteristic in a queue may have an ill effect on another characteristic.

 We shortly discuss   the decay rate of $V$ under some other service disciplines.  Mandjes and Zwart \cite{MandjesZwart} considered the GI/GI/1 \PS\ queue with light-tailed service requests. They show that $\gamma(V_{\PS})$ is equal to $\gamma(L)$, under an additional condition that rules out distributions with bounded support, or extremely light tails. 
Nuyens and Zwart \cite{nuyenszwart} show that $\gamma(V_{\SRPT})$ equals $\gamma(L)$ in the GI/GI/1 \SRPT\ queue, unless there is mass in the endpoint $x_F$ of the service-time distribution. 
Finally, it may be seen that the decay rate $\gamma(V_{\FIFO})$ of the sojourn time in the \FIFO\ system is strictly larger than $\gamma(L)$, see also \cite{MandjesZwart, nuyenszwart}.

For the conditional sojourn time, Nuyens et al.\ \cite{nwz} showed the following result.
\begin{theorem} For all $x$, 
$\gamma(V(x))=\gamma(L_x)$, where $L_x$ is the busy-period length in the queue with generic service time $B\wedge x$.
\end{theorem}

\section{Slowdown}\label{slodo}
We conclude the paper by considering another performance measure for  queueing policies,
 the so-called {\em slowdown}.
The slowdown is a way to measure how fair  jobs are
treated by a service discipline, and is  defined as follows.
\begin{definition}[Slowdown]\label{slowdev} The slowdown $S(x)$ of a job of size $x$ is defined by $S(x)=V(x)/x$. The slowdown $S$ is defined as $S=S(B)$, where $B$ is the generic service time, independent of $S(x)$, and we may write
$ P(S>x)=\int P(S(u)>x)dF(u).$
\end{definition}

Lemma 14 in Nuyens et al.\ \cite{nwz} implies the following result for the asymptotic behaviour of the slowdown.
\begin{theorem}\label{vens} In the \FB\ queue, if $\rho<1$, then
$S(x)$ a.s.\ converges to $1/(1-\rho)$ as $x\to\infty$.
\end{theorem}
Harchol-Balter et al.\ \cite{mor1} proved that in the M/G/1 queue, $\lim_{x\to\infty} S(x)\leq 1/(1-\rho)$ a.s. Hence, the asymptotic slowdown under 
\FB\ is maximal, just as it is under \PS.

The intuitively appealing idea `the larger the service request of a job, the larger his slowdown under \FB\ ' was recently invalidated by Harchol-Balter and Wierman \cite{mor2}. 
\begin{theorem}
The mean slowdown $ES(x)$ is not
monotonically increasing in $x$. In fact, 
$ES(x)$ converges from above to \mbox{$1/(1-\rho)$} as $x\to\infty$.
\end{theorem}
An important  conclusion is   that when one uses the slowdown as a measure of fairness, 
not the longest jobs are treated most unfairly, as is often believed, but certain `medium long' jobs.

Some studies have been done to compare the slowdown in the \FB\ queue with the slowdown in other queues. Let $S_{\FB}$ and $S_{\PS}$ denote the slowdown in two M/G/1 queues with the same arrival rate and service-time distribution.
Theorem 2 of Rai et al.\  \cite{rai2003}  reads
\[ ES_{\FB}\leq \frac{2-\rho}{2-2\rho} ES_{\PS}=\frac{2-\rho}{2(1-\rho)^2}.\]
 Feng and Misra \cite{misra} show that for DFR service-time distributions, the \FB\ discipline minimises the expected slowdown over the class ${\cal D}$.  For more results on the slowdown and (un)fairness we refer to Rai et al.\  \cite{rai2003} and Bansal and Wierman \cite{wierman}.

Rai et al.\  \cite{rai2003} compare, through numerical evaluation,  the slowdown under \FB, \PS\ and \SRPT\  for service times with the so-called {\em high-variability property}. For such service-time distributions less than $1\%$ of the jobs accounts for more than half the load. According to recent studies internet traffic exhibits this property, see Crovella and Bestavros \cite{bestavros}. The numerical study in \cite{rai2003} shows that  for service times with the high-variability property, \FB\ is quite close to \SRPT. Furthermore, a very large percentage of the jobs has a significantly smaller slowdown under \FB\  than under \PS\ or \FIFO.
\bibliography{bibfile}
\appendix
\section{Classes of distributions and stochastic orderings}\label{fiokje}
In this appendix we give an explanation of the  classes of distributions and  stochastic order relations that are used in this paper.

A distribution $F$ with density $f$ belongs to the class DFR (decreasing failure rate), if its failure (or hazard) rate, $f(x)/(1-F(x))$, is decreasing for $x\geq 0$. The class IFR (increasing failure rate) is defined analogously. Alternatively, one may define  the failure rate  only for $x$ in the support of $F$. However, as remarked by  Down and Wu \cite{doug},  this alternative definition is not good enough to prove the optimality for \FB\ in Theorem \ref{dfr} below, and hence it will not be used here.

A density $f$ is called log-convex if  $\log f$ is  convex. By integration, one can show that the class of log-convex densities is a subclass of DFR.
The class of distributions with a log-convex density includes many well-known distributions, for example  Pareto distributions,  gamma distributions with density $f(x)=\lambda^n x^{n-1}\exp(-\lambda x)/\Gamma(n), \lambda\geq 1, x\geq 0$, and Weibull distributions with distribution function
$F(x)=1-\exp(-a x^\beta), \beta\leq 1, x\geq 0$. Densities are called log-concave if $\log f$ is concave. For more results on these classes of distributions  we refer to Shaked and Shanthikumar \cite{shasha}.

Recall that a random variable $X$ is stochastically smaller than $Y$, notation $X\leq_{st} Y$, if $P(X\leq x)\geq P(Y\leq x)$ for all $x$. As a consequence, we can find $X'\stackrel{d}{=}X$ and $Y'\stackrel{d}{=}Y$ such that $P(X'\leq Y')=1$. In fact, this characterisation is equivalent to the definition of $\leq_{st}$.
The stochastic ordering  of  processes is a generalisation of stochastic ordering for random variables, and can be defined similarly, see also Section 4.B.7 of Shaked and Shanthikumar \cite{shasha}: we say that two processes $\{X(t), t\geq 0\}$ and $\{Y(t), t\geq 0\}$
are stochastically ordered, notation $\{X(t), t\geq 0\}\leq_{st} \{Y(t), t\geq 0\}$, if there exist processes $\{\bar{X}(t), t\geq 0\}$ and $\{Y(t), t\geq 0\}$, defined on an common probability space, such that $P(\bar{X}(t)\leq \bar{Y}(t)\  \forall  t)=1$ and 
 \[ \{X(t), t\geq 0\}\stackrel{d}{=}\{\bar{X}(t), t\geq 0\},
 \{Y(t), t\geq 0\}\stackrel{d}{=}\{\bar{Y}(t), t\geq 0\}.\]

\end{document}